\documentclass[a4paper]{article}

\usepackage[margin=1in]{geometry} 

\usepackage[T1]{fontenc}
\newcommand{\changefont}[3]{
\fontfamily{#1} \fontseries{#2} \fontshape{#3} \selectfont}

\changefont{ptm}{m}{n}

\usepackage{setspace} 
\usepackage{graphicx}

\usepackage{amsfonts}
\usepackage{amsmath}
\usepackage{amssymb}
\usepackage{graphics}
\usepackage{mathrsfs}
\usepackage{color}

\newtheorem{theorem}{Theorem}[section]

\newtheorem{lemma}{Lemma}[section]

\newtheorem{definition}{Definition}[section]

\long\def\symbolfootnote[#1]#2{\begingroup%
\def\thefootnote{\fnsymbol{footnote}}\footnote[#1]{#2}\endgroup} 

\begin{document}

\begin{center}
\Large \textbf{Homoclinic and Heteroclinic Trajectories of Differential Equations with Piecewise Constant Arguments of Generalized Type}
\end{center}

\begin{center}
\normalsize \textbf{Mehmet Onur Fen$^{1,}\symbolfootnote[1]{Corresponding Author. E-mail: monur.fen@gmail.com, onur.fen@tedu.edu.tr}$, Fatma Tokmak Fen$^2$} \\
\vspace{0.2cm}
\textit{\textbf{$^1$Department of Mathematics, TED University, 06420 Ankara, Turkey}} \\

\vspace{0.1cm}
\textit{\textbf{$^2$Department of Mathematics, Gazi University, 06560 Ankara, Turkey}} \\
\vspace{0.1cm}
\end{center}

\vspace{0.3cm}

\begin{center}
\textbf{Abstract} 
\end{center}

\noindent Quasilinear systems with piecewise constant arguments of generalized type are under investigation from the asymptotic point of view. The systems have discontinuous right-hand sides which are identified via a discrete-time map. It is rigorously proved that  homoclinic and heteroclinic solutions are generated, and they are taken into account in the functional sense. The Banach fixed point theorem is used for the verification. The hyperbolic set of solutions is also discussed, and an example supporting the theoretical findings is provided.
\vspace{-0.2cm}

\noindent\ignorespaces

\vspace{0.3cm}
 
\noindent\ignorespaces \textbf{Keywords:} piecewise constant argument, homoclinic solution, heteroclinic solution, hyperbolic set of functions 

\vspace{0.2cm}

\noindent\ignorespaces \textbf{2020 Mathematics Subject Classification:} 34K16, 34K34, 37C29

\vspace{0.6cm}


\section{Introduction} \label{sec1}

In the present study, we take into account differential equations with piecewise constant argument of generalized type (EPCAG) of the form
\begin{eqnarray} \label{mainsystem}
z'(t) = Az(t) + f(t,z(t), z(\gamma (t))) + g(t,\alpha),
\end{eqnarray}
where $t\in\mathbb R$, $A \in \mathbb R^{m\times m}$ is a matrix whose all eigenvalues have negative real parts, $f:\mathbb R \times \mathbb R^m \times \mathbb R^m \to \mathbb R^m$ is a function which is continuous in all of its arguments, $\gamma(t) = \zeta_k$ and $g(t,\alpha) = \alpha_k$ for $\theta_k \leq t < \theta_{k+1}$, $k\in\mathbb Z$. Here, $\{\theta_k\}_{k\in \mathbb Z}$ and $\{\zeta_k\}_{k\in \mathbb Z}$ are fixed sequences of real numbers such that $\theta_k \leq \zeta_k\leq\theta_{k+1}$ and $\theta_{k+1}-\theta_k=\omega$, $k\in\mathbb Z$, for some $\omega>0$. Additionally, $\alpha=\{\alpha_k\}_{k\in\mathbb Z}$ is an orbit of the  discrete-time equation
\begin{eqnarray} \label{discretemap}
\eta_{k+1} = F(\eta_{k}),
\end{eqnarray}
in which $F:\Gamma \to \Gamma$ is a continuous function and $\Gamma$ is a bounded subset of $\mathbb R^m$.

Researches on EPCAG were initiated by Akhmet \cite	{Akhmet2007}-\cite{Akhmet2008c}, and this type of differential equations is a general case of the ones discussed in the studies \cite{Cooke84,Wiener88,Wienerbook}, where the greatest integer function is used for the construction of the argument. Because both continuous and discrete arguments take place in EPCAG, such equations belong to a class of hybrid systems \cite{Alwan13}. The reader is referred to the studies \cite{Cui24,Karakoc17,Ozturk12,Qiuxiang06,Yu2016} for several applications of systems with piecewise constant arguments.

The main purpose of the present paper is the verification of homoclinic and heteroclinic solutions of system (\ref{mainsystem}). To the best of our knowledge, this is the first time in the literature that the existence of such solutions is reported for differential equations with piecewise constant arguments, and especially for EPCAG. Homoclinic and heteroclinic trajectories play an important role for the verification of chaos \cite{Silva93} and are observable in real world phenomena such as mechanical systems, chemical reactions, economics, and predator-prey models \cite{Agliari11,Petersen88,Mikhlin03,Wang19}. We would like to point out that systems of the form (\ref{mainsystem}) comprise a large class of differential equations including the ones with alternately of retarded and advanced types. For instance, by taking $\theta_k=\zeta_k=k$, $k\in\mathbb Z$, in  (\ref{mainsystem})  the equation 
\begin{eqnarray} \label{particular1}
	z'(t) = Az(t) + f(t,z(t), z([t])) + g(t,\alpha),
\end{eqnarray}
is obtained, whereas the choices $\theta_k=m_1k-m_2$, $\zeta_k=m_1k$, $k\in\mathbb Z$, make it possible to attain 
\begin{eqnarray} \label{particular2}
	z'(t) = Az(t) + f \left( t,z(t), z\left( m_1 \left[ \displaystyle \frac{t+m_2}{m_1} \right] \right)  \right) + g(t,\alpha),
\end{eqnarray}
where $[.]$ stands for the greatest integer function and $m_1,m_2$ are positive integers with $m_1>m_2$. For that reason, the theoretical results obtained in the present paper are new also for systems of the form (\ref{particular1}) and (\ref{particular2}).

Based on the Melnikov method and topological horseshoes, Dong and Li \cite{Dong25} demonstrated the existence of homoclinic and heteroclinic solutions for systems with small perturbations. Piecewise linear systems possessing these types of solutions were investigated in paper \cite{Medrano03} numerically. The Fishing principle is another technique that is useful for their detection \cite{Leonov14}. On the other hand, Melnikov type conditions were derived in \cite{Battelli24} to deduce the persistence of heteroclinic solutions in singularly perturbed discontinuous models. The key difference of our study is the presence of piecewise constant arguments in the model. The systems discussed in \cite{Battelli24,Dong25,Leonov14,Medrano03} do not comprise piecewise constant arguments. Our proof technique is also distinctive such that we make benefit of the uniformly bounded solutions generated by (\ref{mainsystem}) as well as the Banach fixed point theorem. Moreover, we utilize the definitions of homoclinic and heteroclinic motions in functional sense. 

The rest of the paper is organized in the following way. In Section \ref{sec2}, the required assumptions on system (\ref{mainsystem}) and the definitions of homoclinic and heteroclinic motions are provided. The bounded solutions are also discussed in that section. The main results of the present research are given in Section \ref{sec3}, where we investigate the stable and unstable sets of bounded solutions as well as the existence of homoclinic and heteroclinic motions. Section \ref{sec4}, on the other hand, is devoted to an example. The logistic map has a rich dynamical structure including homoclinic and heteroclinic orbits \cite{Avrutin15}, and it is utilized for the construction of the example. Finally, concluding remarks with a possible future problem are mentioned in Section \ref{sec5}.

 \section{Preliminaries} \label{sec2}
 
Let us denote by $\Lambda$ the set of all orbits $\alpha=\{ \alpha_k \}_{k\in\mathbb Z}$ of the discrete-time equation (\ref{discretemap}) which take place in the bounded set $\Gamma$. The description for a solution of (\ref{mainsystem}) defined on the whole real axis is as follows.
 
\begin{definition} (\cite{Akhmet2008a})
For a fixed sequence $\alpha = \{ \alpha_k \}_{k\in\mathbb Z} \in \Lambda$, a continuous function $\psi:\mathbb R \to \mathbb R^m$ is a solution of system (\ref{mainsystem}) if:
\begin{itemize}
\item[\textbf{(i)}] the derivative $\psi'(t)$ exists at each $t\in\mathbb R$ with the possible exception of the points $\theta_k$, $k\in\mathbb Z$, where the one-sided derivatives exist;
\item[\textbf{(ii)}] $\psi(t)$ satisfies the equation (\ref{mainsystem}) on each interval $(\theta_k,\theta_{k+1})$, $k\in\mathbb Z$, and it holds for the right derivative of $\psi(t)$ at the points $\theta_k$, $k\in\mathbb Z$.
 \end{itemize}
 \end{definition}
 
In what follows the usual Euclidean norm for vectors and the spectral norm for square matrices are utilized.

 Because the matrix $A$ in system (\ref{mainsystem}) admits eigenvalues with negative real parts, there exist a number $N\geq 1$ and a number $\lambda >0$ such that $\left\|e^{At} \right\|  \leq N e^{-\lambda t}$ for all $t\geq 0$.
 
 The following assumptions on system (\ref{mainsystem}) are required.
 \begin{itemize}
\item[\textbf{(A1)}] $\displaystyle \sup_{(t,x,y) \in \mathbb R \times \mathbb R^m \times \mathbb R^m } \left\|f(t,x,y) \right\|\leq M_f $ for some   positive number $M_f$;
\item[\textbf{(A2)}] There is a positive number $L_1$ such that $\left\|f(t,x_1,y) - f(t,x_2,y) \right\| \leq L_1 \left\|x_1-x_2 \right\|$ for every $t\in\mathbb R$, $x_1,x_2,y\in\mathbb R^m$; 
\item[\textbf{(A3)}] There is a positive number $L_2$ such that $\left\|f(t,x,y_1) - f(t,x,y_2) \right\| \leq L_2 \left\|y_1-y_2 \right\|$ for every $t\in\mathbb R$, $x,y_1,y_2\in\mathbb R^m$;
\item[\textbf{(A4)}] $N(L_1+L_2) <\lambda$.
\end{itemize}

One can confirm using the techniques presented in the studies \cite{Akhmet2007,Hartman64} that for a fixed orbit $\alpha=\{\alpha_k\}_{k\in\mathbb Z}$ of (\ref{discretemap}) a bounded function $z:\mathbb R\to \mathbb R^m$ is a solution of  system (\ref{mainsystem}) if and only if the integral equation
$$
z(t) = \displaystyle \int_{-\infty}^t e^{A(t-s)} \left[ f(s,z(s), z(\gamma (s))) + g(s,\alpha)\right] ds
$$
is fulfilled.
In the remaining parts of the paper we denote 
\begin{eqnarray*}
M_F = \displaystyle \sup_{s \in \Omega} \left\|F(s) \right\|
\end{eqnarray*}
and 
\begin{eqnarray} \label{numbermphi}
M_{\phi} = \displaystyle \frac{N(M_f +M_F)}{\lambda}.
\end{eqnarray}

The subsequent assertion deals with the bounded solutions of system (\ref{mainsystem}). 

\begin{lemma} \label{bddsolnlemma}
Suppose that the assumptions $(A1)-(A4)$ hold. Then, for every sequence $\alpha=\{\alpha_k\}_{k\in\mathbb Z} \in \Lambda$, system (\ref{mainsystem}) possesses a unique solution $\phi_{\alpha}(t)$ which is bounded on the real axis such that $\displaystyle \sup_{t\in\mathbb R}\left\| \phi_{\alpha}(t)\right\| \leq  M_{\phi}$.
\end{lemma}

Lemma \ref{bddsolnlemma} can be proved making use of the Banach fixed point theorem, and for the convenience of the reader the verification is provided in the Appendix. 

For a fixed sequence $\alpha=\{ \alpha_k \}_{k\in\mathbb Z} \in \Lambda$, the sets
$$W^s(\alpha) = \left\lbrace \beta \in \Lambda : \, \left\| \alpha_k - \beta_k \right\| \to 0 \,\textrm{as}\, k\to\infty    \right\rbrace $$
and
$$W^u(\alpha) = \left\lbrace \beta \in \Lambda : \, \left\| \alpha_k - \beta_k \right\| \to 0 \,\textrm{as}\, k\to -\infty    \right\rbrace $$
are respectively called the stable and unstable sets of $\alpha$ \cite{Akhmet2008b,Akhmet2010}.

\begin{definition} (\cite{Akhmet2008b,Akhmet2010}) \label{defn1}
A sequence $\beta \in \Lambda$ is called homoclinic to a sequence $\alpha\in\Lambda$ if $\alpha\neq\beta$ and $\beta \in W^s(\alpha) \cap W^u(\alpha)$. Moreover, $\beta \in \Lambda$ is called heteroclinic to sequences $\alpha^1,\alpha^2\in\Lambda$ if $\alpha^1\neq \alpha^2$, $\alpha^1\neq\beta$, $\alpha^2\neq\beta$,  and $\beta \in W^s(\alpha^1) \cap W^u(\alpha^2)$. 
\end{definition}

\begin{definition}  (\cite{Akhmet2008b,Akhmet2010})  \label{defn2}
	The set $\Lambda$ is said to be hyperbolic if for each $\alpha \in \Lambda$ the stable and unstable sets of $\alpha$ contain at least one element different from $\alpha$. 
\end{definition}

In compliance with the result of Lemma \ref{bddsolnlemma}, let us denote by $\Omega$ the set of all bounded solutions of (\ref{mainsystem}), i.e.,
$\Omega = \left\lbrace \phi_{\alpha}(t) : \alpha \in \Lambda \right\rbrace$. It is worth noting that two elements $\phi_{\alpha}(t), \phi_{\beta}(t)$ of $\Omega$ are equal  if and only if  $\alpha_k= \beta_k$ for each $k\in\mathbb Z$. If $\phi_{\alpha}(t)$ is a function in $\Omega$, then its stable set is defined by
$$W^s\left( \phi_{\alpha}(t)\right) = \left\lbrace  \phi_{\beta}(t)\in\Omega: \left\|\phi_{\alpha}(t) -\phi_{\beta}(t)\right\| \to 0 \textrm{ as } t\to\infty \right\rbrace,$$
whereas its unstable set is given by 
$$W^u\left( \phi_{\alpha}(t)\right) = \left\lbrace  \phi_{\beta}(t)\in\Omega: \left\|\phi_{\alpha}(t) -\phi_{\beta}(t)\right\| \to 0 \textrm{ as } t\to -\infty \right\rbrace.$$

\begin{definition} (\cite{Akhmet2008b,Akhmet2010}) \label{defn3}
 A function $\phi_{\beta}(t) \in \Omega$ is called homoclinic to a function $\phi_{\alpha}(t) \in \Omega$ if $\alpha\neq\beta$ and $\phi_{\beta}(t) \in W^s(\phi_{\alpha}(t)) \cap W^u(\phi_{\alpha}(t))$. Moreover, $\phi_{\beta}(t)\in \Omega$ is called heteroclinic to functions $\phi_{\alpha^1}(t)$, $\phi_{\alpha^2}(t)\in\Omega$ if $\alpha^1\neq\alpha^2$, $\alpha^1\neq\beta$, $\alpha^2\neq\beta$, and $\phi_{\beta}(t) \in W^s(\phi_{\alpha^1}(t)) \cap W^u(\phi_{\alpha^2}(t))$. 
\end{definition}

\begin{definition} (\cite{Akhmet2008b,Akhmet2010}) \label{defn4}
The set $\Omega$ is called hyperbolic if for every $\alpha \in \Lambda$ the sets $W^s(\phi_{\alpha}(t))$ and $W^u(\phi_{\alpha}(t))$ contain at least one element different from $\phi_{\alpha}(t)$. 
\end{definition}

Definitions \ref{defn3} and \ref{defn4} are utilized in the next section to achieve the main results of the present study.

\section{Main Results} \label{sec3}

To approve the existence of homoclinic and heteroclinic solutions of system (\ref{mainsystem}), the following assumption is needed.
 \begin{itemize}
\item[\textbf{(A5)}] $\displaystyle \frac{N}{\lambda} \left(  2   L_1 + \frac{ L_2 e^{\lambda \omega/2} \left( e^{\lambda \omega}-1\right) }{ 1-e^{-\lambda \omega/2} }\right) <1$.
\end{itemize}
It is worth noting that the assumption $(A5)$ implies $(A4)$.

A result concerning the stable sets for bounded solutions of (\ref{mainsystem}) is given in the next assertion.
\begin{lemma} \label{lemmahomoclinic1}
Suppose that the assumptions $(A1)-(A3)$ and $(A5)$ are fulfilled. If $\alpha=\left\lbrace \alpha_k\right\rbrace_{k\in\mathbb Z}$ and $\beta=\left\lbrace \beta_k\right\rbrace_{k\in\mathbb Z}$ are elements of $\Lambda$ such that $\beta \in W^s(\alpha)$, then $\phi_{\beta}(t) \in W^s (\phi_{\alpha} (t))$. 
\end{lemma}

\noindent \textbf{Proof.} 
Let us define the numbers $$R_1 = 2 N M_{\phi} \left(1-\frac{2 N L_1}{\lambda} - \frac{N L_2 e^{\lambda \omega/2} \left( e^{\lambda \omega}-1\right) }{\lambda \left( 1-e^{-\lambda \omega/2}\right) } \right)^{-1},$$
$$R_2 = \frac{N}{\lambda}  \left( 1- \frac{NL_1}{\lambda}  - \frac{N L_2 e^{\lambda \omega}}{\lambda} \right)^{-1},$$ and suppose that $\sigma$ is a positive number with
$$\sigma < \displaystyle \frac{1}{R_1 +R_2}.$$ The numbers $R_1$ and $R_2$ are positive in accordance with the assumption $(A5)$.

Fix a number $\varepsilon>0$. Because the sequence $\beta=\left\lbrace \beta_k\right\rbrace_{k\in\mathbb Z}$ belongs to $W^s(\alpha)$, there exists an integer $k_0$ such that $$\left\| \alpha_k -\beta_k\right\| < \sigma \varepsilon, \, k\geq k_0.$$  

One can confirm that the function $u(t) = \phi_{\alpha}(t) - \phi_{\beta}(t)$ is a solution of the system
\begin{eqnarray} \label{wsystem}
u'(t) = Au(t) + f\left( t, u(t) + \phi_{\beta}(t), u(\gamma(t)) + \phi_{\beta}(\gamma(t))  \right) -  f\left( t,  \phi_{\beta}(t),   \phi_{\beta}(\gamma(t))  \right) + g(t,\alpha) - g(t,\beta)
\end{eqnarray}
and satisfies the integral equation
\begin{eqnarray*}
u(t) & =& e^{A \left( t-\theta_{k_0} \right) } \left( \phi_{\alpha}(\theta_{k_0})  - \phi_{\beta}(\theta_{k_0}) \right) 	\\
& +& \displaystyle \int_{\theta_{k_0}}^t e^{A(t-s)} \big[ f\left( s, u(s) + \phi_{\beta}(s), u(\gamma(s)) + \phi_{\beta}(\gamma(s))  \right) 
-  f\left( s,  \phi_{\beta}(s),   \phi_{\beta}(\gamma(s))  \right)  \big] ds \\
&+&  \displaystyle \int_{\theta_{k_0}}^t e^{A(t-s)} \left[g(s,\alpha)  - g(s,\beta) \right] ds. 
\end{eqnarray*}

Let us denote by $\mathcal D$ the set of all uniformly bounded continuous functions $u:\mathbb R \to \mathbb R^m$ satisfying 
$$\left\| u(t)\right\| \leq R_1 e^{-\lambda (t-\theta_{k_0})/2} +R_2 \sigma \varepsilon, \  t\geq \theta_{k_0},$$
 and $\left\| u\right\|_{\infty} \leq H$, where $\left\| u\right\|_{\infty} =\displaystyle \sup_{t\in\mathbb R} \left\| u(t)\right\|$ and $H= 2N\displaystyle \left(  M_{\phi} + \frac{M_f + M_F}{\lambda}\right).$

We take into account the operator $T$ defined on the set $\mathcal D$ by means of the equation
$$
(T(u))(t) = \begin{cases}
	\phi_{\alpha}(t) - \phi_{\beta}(t),  ~ \textrm{if }  t < \theta_{k_0}, \\
	e^{A \left( t-\theta_{k_0} \right) } \left( \phi_{\alpha}(\theta_{k_0})  - \phi_{\beta}(\theta_{k_0}) \right) 	\\
	  +  \displaystyle \int_{\theta_{k_0}}^t e^{A(t-s)} \big[ f\left( s, u(s) + \phi_{\beta}(s), u(\gamma(s)) + \phi_{\beta}(\gamma(s))  \right) 
	-  f\left( s,  \phi_{\beta}(s),   \phi_{\beta}(\gamma(s))  \right)  \big] ds \\
	 +   \displaystyle \int_{\theta_{k_0}}^t e^{A(t-s)} \left[g(s,\alpha) - g(s,\beta) \right] ds , ~  \textrm{if }  t \geq \theta_{k_0}.
\end{cases}
$$

Firstly, we will demonstrate that $T(\mathcal D) \subseteq \mathcal D$. Suppose that $u(t)$ is an element of $\mathcal D$, and fix an arbitrary  $t\geq \theta_{k_0}$. Then, we have
\begin{eqnarray*}
\left\|(T(u))(t) \right\| & \leq & Ne^{-\lambda \left( t-\theta_{k_0}\right) } \left\| \phi_{\alpha} \left(\theta_{k_0} \right) - \phi_{\beta} \left(\theta_{k_0} \right) \right\| \\
&+& \displaystyle \int_{\theta_{k_0}}^t Ne^{-\lambda(t-s)} \big\|  f\left(s, u(s) + \phi_{\beta}(s), u(\gamma(s)) + \phi_{\beta}(\gamma(s))  \right) 
\\ && - f\left(s,\phi_{\beta}(s), u(\gamma(s)) + \phi_{\beta}(\gamma(s)) \right)  \big\| ds \\
&+& \displaystyle \int_{\theta_{k_0}}^t Ne^{-\lambda(t-s)}  \big\|  f\left(s, \phi_{\beta}(s), u(\gamma(s)) + \phi_{\beta}(\gamma(s)) \right) 
 -  f\left( s, \phi_{\beta}(s), \phi_{\beta}(\gamma(s))  \right)  \big\| ds \\
&+& \displaystyle \int_{\theta_{k_0}}^t Ne^{-\lambda(t-s)} \left\|g(s,\alpha) - g(s,\beta) \right\| ds \\
&\leq &  2N M_{\phi} e^{-\lambda \left( t-\theta_{k_0}\right)} 
+ N L_1 \displaystyle \int_{\theta_{k_0}}^t e^{-\lambda(t-s)} \left\|u(s) \right\|ds  \\
&+& N L_2 \displaystyle \int_{\theta_{k_0}}^t e^{-\lambda(t-s)} \left\|u(\gamma(s)) \right\|ds  + N \sigma \varepsilon \displaystyle \int_{\theta_{k_0}}^t  e^{-\lambda(t-s)} ds \\
&\leq& 2N M_{\phi} e^{-\lambda ( t-\theta_{k_0})/2} 
+ \displaystyle \frac{N \sigma \varepsilon}{\lambda} \left( 1- e^{-\lambda (t-\theta_{k_0})} \right)  \\
&+& N L_1 \displaystyle \int_{\theta_{k_0}}^t e^{-\lambda(t-s)} \left\|u(s) \right\|ds 
+N L_2 \displaystyle \int_{\theta_{k_0}}^t e^{-\lambda(t-s)} \left\|u(\gamma(s)) \right\|ds. 
\end{eqnarray*}
One can confirm that
\begin{eqnarray} \label{proofineq1}
\displaystyle \int_{\theta_{k_0}}^t e^{-\lambda(t-s)} \left\|u(s) \right\|ds 
& \leq & \displaystyle \int_{\theta_{k_0}}^t e^{-\lambda(t-s)}  \left( R_1 e^{-\lambda (s-\theta_{k_0})/2} +R_2 \sigma \varepsilon\right)  ds \nonumber \\
& = & \displaystyle \frac{2 R_1}{\lambda} e^{-\lambda  (t-\theta_{k_0}  )/2 } \left(1-e^{-\lambda (t-\theta_{k_0})/2 } \right) + \frac{R_2\sigma \varepsilon}{\lambda} \left( 1-e^{-\lambda  (t-\theta_{k_0}  ) }\right)  \nonumber \\
& < & \displaystyle \frac{2 R_1}{\lambda} e^{-\lambda  (t-\theta_{k_0}  )/2 }  + \frac{R_2 \sigma \varepsilon}{\lambda}.  
\end{eqnarray}
Next, let $p$ be the integer such that $\theta_{k_{0}+p} < t \leq \theta_{k_{0} + p + 1}$.  In this case we have
\begin{eqnarray} \label{proofineqp}
	\displaystyle \int_{\theta_{k_0}}^t  e^{-\lambda(t-s)} \left\|u(\gamma(s)) \right\|ds \leq \displaystyle \sum_{j=0}^{p} \left( \int_{\theta_{k_{0}+j}}^{\theta_{k_{0}+j+1}}  e^{-\lambda(t-s)} \left\|u(\zeta_{k_0+j}) \right\|ds \right). 
\end{eqnarray}
It can be verified that
\begin{eqnarray*}
\int_{\theta_{k_{0}+j}}^{\theta_{k_{0}+j+1}}  e^{-\lambda(t-s)} \left\|u(\zeta_{k_0+j}) \right\|ds 
&\leq&  \int_{\theta_{k_{0}+j}}^{\theta_{k_{0}+j+1}}  e^{-\lambda(t-s)}  \left( R_1 e^{-\lambda ( \zeta_{k_0+j}-\theta_{k_0})/2} +R_2 \sigma \varepsilon\right)  ds \\
&=& \displaystyle \frac{e^{\lambda \omega}-1}{\lambda} \left( R_1 e^{-\lambda ( \zeta_{k_0+j}-\theta_{k_0})/2} +R_2 \sigma \varepsilon\right)  e^{-\lambda (t-\theta_{k_0+j})} \\
&\leq& \displaystyle \frac{e^{\lambda \omega}-1}{\lambda} \left(R_1 e^{-\lambda \omega j/2} + R_2 \sigma \varepsilon \right) e^{-\lambda (p-j) \omega}  \\
&=& \displaystyle \frac{\left( e^{\lambda \omega} -1\right) e^{-\lambda p \omega} }{\lambda} \left(R_1 e^{\lambda \omega j /2} + R_2 \sigma \varepsilon e^{\lambda \omega j} \right). 
\end{eqnarray*}
We attain using (\ref{proofineqp}) together with the inequality $p\geq (t-\theta_{k_0})/\omega-1$ that 
\begin{eqnarray} \label{proofineq2}
\displaystyle \int_{\theta_{k_0}}^t  e^{-\lambda(t-s)} \left\|u(\gamma(s)) \right\|ds 
&\leq& \displaystyle \frac{\left( e^{\lambda \omega} -1\right) e^{-\lambda p \omega} }{\lambda} \left(R_1 \sum_{j=0}^p e^{\lambda \omega j /2} + R_2 \sigma \varepsilon \sum_{j=0}^p e^{\lambda \omega j} \right) \nonumber \\
&=& \displaystyle \frac{\left( e^{\lambda \omega} -1\right) e^{-\lambda p \omega} }{\lambda}  \left(R_1 \frac{e^{\lambda  (p+1)\omega /2} -1}{e^{\lambda \omega/2}-1} + R_2 \sigma \varepsilon \frac{e^{\lambda  (p+1)\omega}-1}{e^{\lambda \omega} -1} \right) \nonumber \\
&<& \displaystyle \frac{R_1 e^{\lambda \omega/2} \left(e^{\lambda\omega} -1 \right)}{\lambda \left( 1- e^{-\lambda \omega/2}\right)}  e^{-\lambda (t-\theta_{k_0})/2} + \frac{R_2 e^{\lambda \omega} \sigma \varepsilon}{\lambda}.
\end{eqnarray}
The inequalities (\ref{proofineq1}) and (\ref{proofineq2}) yield for every $t \geq \theta_{k_0}$ that
\begin{eqnarray*}
\left\|(T(u))(t) \right\| &<& 2NM_{\phi} e^{-\lambda (t-\theta_{k_0}) /2} + \displaystyle \frac{N\sigma \varepsilon}{\lambda}  
 +  \displaystyle \frac{2 R_1 N L_1}{\lambda} e^{-\lambda (t-\theta_{k_0})/2} + \frac{R_2 N L_1 \sigma \varepsilon}{\lambda} \\
&+& \displaystyle \frac{R_1 N L_2 e^{\lambda \omega/2} \left(e^{\lambda \omega -1} \right) }{\lambda \left( 1- e^{-\lambda \omega/2}\right) } + \frac{R_2 N L_2 e^{\lambda \omega} \sigma \varepsilon}{\lambda} \\
&=& R_1 e^{-\lambda (t-\theta_{k_0})/2} + R_2 \sigma \varepsilon.
\end{eqnarray*}

Now, we will discuss the boundedness of $T(u)$. The inequality $$\left\| (T(u)) (t)\right\| = \left\| \phi_{\alpha}(t) - \phi_{\beta} (t)\right\| \leq 2M_{\phi}\leq H$$ holds whenever $t<\theta_{k_0}$.
On the other hand, if $t \geq \theta_{k_0}$, then we have
\begin{eqnarray*}
\left\| (T(u)) (t)\right\| &\leq& N e^{-\lambda (t-\theta_{k_0})} \left\| \phi_{\alpha} (\theta_{k_0}) -\phi_{\beta} (\theta_{k_0}) \right\|  \\
& +& \displaystyle \int_{\theta_{k_0}}^t N e^{-\lambda(t-s)}  \left\|  f\left( s, u(s) + \phi_{\beta}(s), u(\gamma(s)) + \phi_{\beta}(\gamma(s))  \right) 
-  f\left( s,  \phi_{\beta}(s),   \phi_{\beta}(\gamma(s))  \right)   \right\|  ds \\
&+&  \displaystyle \int_{\theta_{k_0}}^t N e^{-\lambda(t-s)}  \left\| g(s,\alpha)  - g(s,\beta)  \right\|  ds \\
&\leq & 2NM_{\phi} e^{-\lambda (t-\theta_{k_0})} + \frac{2N(M_f +M_F)}{\lambda} \left( 1-e^{-\lambda (t-\theta_{k_0})}\right) \\
&\leq& 2NM_{\phi} + \frac{2N(M_f +M_F)}{\lambda} =H. 
\end{eqnarray*}
Accordingly, we have $\left\| T(u)\right\|_{\infty}  \leq H$. Thus, $T(\mathcal D) \subseteq \mathcal D$.

To approve that $T$ is a contraction mapping, let us take two elements $u_1(t)$ and $u_2(t)$ of $\mathcal D$. For $t<\theta_{k_0}$ we have $(T(u_1)) (t) - (T(u_2)) (t)=0$. If $t \geq \theta_{k_0}$, then using the equation
\begin{eqnarray*}
(T(u_1)) (t) - (T(u_2)) (t) &=& \displaystyle \int_{\theta_{k_0}}^t e^{A(t-s)} \big[ f\left( s, u_1(s) + \phi_{\beta}(s), u_1(\gamma(s)) + \phi_{\beta}(\gamma(s))  \right) \\
&& -  f\left( s, u_2(s) + \phi_{\beta}(s), u_2(\gamma(s)) + \phi_{\beta}(\gamma(s))  \right)\big] ds
\end{eqnarray*}
we obtain that
\begin{eqnarray*}
\left\| (T(u_1)) (t) - (T(u_2)) (t) \right\| 
&\leq & \displaystyle \int_{\theta_{k_0}}^t N L_1 e^{-\lambda (t-s)} \left\|u_1(s) - u_2(s) \right\| ds \\
&+&  \displaystyle \int_{\theta_{k_0}}^t N L_2 e^{-\lambda (t-s)} \left\|u_1(\gamma (s)) - u_2(\gamma (s)) \right\| ds \\
&\leq& \displaystyle \frac{N(L_1+L_2)}{\lambda} \left( 1-e^{-\lambda (t-\theta_{k_0})}\right) \left\| u_1-u_2\right\|_{\infty}. 
\end{eqnarray*}
Thus, $\left\|T(u_1) -T(u_2) \right\|_{\infty} \leq \displaystyle \frac{N(L_1+L_2)}{\lambda} \left\|u_1-u_2 \right\|_{\infty}$, and accordingly, $T$ is a contraction mapping.

Owing to the uniqueness of solutions for system (\ref{wsystem}), the function $u(t)=\phi_{\alpha} (t) - \phi_{\beta} (t)$ is the unique fixed point of the operator $T$. Let us denote
$$
u_0(t) = \begin{cases}
	\phi_{\alpha}(t) - \phi_{\beta}(t),  ~ \textrm{if }  t < \theta_{k_0}, \\
	e^{A \left( t-\theta_{k_0} \right) } \left( \phi_{\alpha}(\theta_{k_0})  - \phi_{\beta}(\theta_{k_0}) \right), ~  \textrm{if }  t \geq \theta_{k_0}.
\end{cases}
$$
The function $u_0(t)$ belongs to $\mathcal D$, and the sequence of functions $\{u_k (t)\}$ with $u_{k+1} = T(u_k)$, $k=0,1,2,\ldots$, converges on the real axis to  $u(t)=\phi_{\alpha} (t) - \phi_{\beta} (t)$. Accordingly, the inequality
\begin{eqnarray*}  
\left\| \phi_{\alpha}(t) -\phi_{\beta}(t)\right\|  \leq  R_1 e^{-\lambda (t-\theta_{k_0})} + R_2 \sigma \varepsilon  
\end{eqnarray*} 
is satisfied for $t\geq \theta_{k_0}$. If $t$ is sufficiently large such that $$t \geq \displaystyle \max \left\lbrace \theta_{k_0} , \theta_{k_0} + \frac{2}{\lambda} \ln \left(\frac{1}{\sigma \varepsilon} \right)  \right\rbrace,$$ then
\begin{eqnarray*}
\left\| \phi_{\alpha}(t) -\phi_{\beta}(t)\right\|  
 \leq   (R_1 +R_2) \sigma \varepsilon < \varepsilon.
\end{eqnarray*}
Hence, $\left\| \phi_{\alpha}(t) -\phi_{\beta}(t)\right\| \to 0$ as $t\to \infty$. Consequently, $\phi_{\beta}(t) \in W^s (\phi_{\alpha} (t))$.  $\square$

The next lemma is devoted to the unstable sets of functions in $\Omega$.

\begin{lemma} \label{lemmahomoclinic2}
Suppose that the assumptions $(A1)-(A4)$ are fulfilled. If $\alpha=\left\lbrace \alpha_k\right\rbrace_{k\in\mathbb Z}$ and $\beta=\left\lbrace \beta_k\right\rbrace_{k\in\mathbb Z}$ are elements of $\Lambda$ such that $\beta \in W^u(\alpha)$, then $\phi_{\beta}(t) \in W^u(\phi_{\alpha}(t))$. 
\end{lemma}

\noindent \textbf{Proof.} Let a positive number $\varepsilon$ be given, and suppose that $\eta$ is a real number such that $$0 < \eta < \displaystyle \frac{\lambda - N(L_1 +L_2)}{N}.$$ 
Since $\beta \in W^u(\alpha)$, there exists an integer $k_0$ such that 
\begin{eqnarray} \label{secondproof1}
\left\| \alpha_k-\beta_k\right\| <\eta \varepsilon
\end{eqnarray}
whenever $k<k_0$.

For every $t \leq \theta_{k_0}$, making use of the equations
$$\phi_{\alpha}(t) = \displaystyle \int_{-\infty}^{t} e^{A(t-s)} \left[ f(s,\phi_{\alpha}(s), \phi_{\alpha} (\gamma(s))) + g(s,\alpha)\right] ds $$
and
$$\phi_{\beta}(t) = \displaystyle \int_{-\infty}^{t} e^{A(t-s)} \left[ f(s,\phi_{\beta}(s), \phi_{\beta}(\gamma(s))) + g(s,\alpha)\right] ds, $$
we attain that
\begin{eqnarray*}
\left\| \phi_{\alpha}(t)-\phi_{\beta}(t)\right\|  
&\leq& N L_1 \displaystyle \int_{-\infty}^t  e^{-\lambda (t-s)} \left\|\phi_{\alpha}(s)  -\phi_{\beta}(s) \right\|ds \\
&+& N L_2 \displaystyle \int_{-\infty}^t  e^{-\lambda (t-s)} \left\|\phi_{\alpha}(\gamma(s))  -\phi_{\beta}(\gamma(s)) \right\|ds \\
&+& N \displaystyle \int_{-\infty}^t e^{-\lambda (t-s)} \left\|g(s,\alpha) -g(s,\beta) \right\|ds. 
\end{eqnarray*}

Fix an arbitrary $t\leq \theta_{k_0}$, and let $q$ be the integer with $q \leq k_0$ such that $\theta_{q-1} < t \leq \theta_q$. Then, we have
\begin{eqnarray*}
&& \displaystyle \int_{-\infty}^t  e^{-\lambda (t-s)} \left\|\phi_{\alpha}(\gamma(s))  -\phi_{\beta}(\gamma(s)) \right\|ds \\
&& = \int_{\theta_{q-1}}^t e^{-\lambda (t-s)} \left\| \phi_{\alpha} (\zeta_{q-1}) - \phi_{\beta} (\zeta_{q-1})  \right\| ds 
+ \displaystyle \sum_{j=1}^{\infty} \left( \int_{\theta_{q-j-1}}^{\theta_{q-j}} e^{-\lambda (t-s)}  \left\| \phi_{\alpha}(\zeta_{q-j-1}) -\phi_{\beta}(\zeta_{q-j-1})  \right\|  ds \right)  \\
&& \leq  \displaystyle \sup_{t \leq \theta_{k_0}} \left\|\phi_{\alpha}(t) -\phi_{\beta}(t) \right\| \left(\int_{\theta_{q-1}}^t e^{-\lambda (t-s)}ds + \sum_{j=1}^{\infty} \left( \int_{\theta_{q-j-1}}^{\theta_{q-j}} e^{-\lambda (t-s)} ds  \right)  \right)  \\
&& =  \displaystyle \frac{1}{\lambda} \sup_{t \leq \theta_{k_0}} \left\|\phi_{\alpha}(t) -\phi_{\beta}(t) \right\|.
\end{eqnarray*}
Moreover, one can confirm by means of (\ref{secondproof1}) that
\begin{eqnarray*}
&& \displaystyle \int_{-\infty}^t e^{-\lambda (t-s)} \left\|g(s,\alpha) -g(s,\beta) \right\|ds \\
&& = \displaystyle \int_{\theta_{q-1}}^t e^{-\lambda (t-s)} \left\| \alpha_{q-1} - \beta_{q-1}\right\| ds  
 + \displaystyle \sum_{j=1}^{\infty} \left( \int_{\theta_{q-j-1}}^{\theta_{q-j}}  e^{-\lambda(t-s)} \left\|\alpha_{q-j-1} -\beta_{q-j-1} \right\| ds \right) \\
 && < \displaystyle \int_{\theta_{q-1}}^t \eta \varepsilon e^{-\lambda (t-s)}  ds  
 + \displaystyle \sum_{j=1}^{\infty} \left( \int_{\theta_{q-j-1}}^{\theta_{q-j}} \eta \varepsilon e^{-\lambda(t-s)} ds \right) \\
 &&=\frac{\eta \varepsilon}{\lambda}. 
\end{eqnarray*}
Now, utilizing the estimation
\begin{eqnarray*}
\displaystyle \int_{-\infty}^t   e^{-\lambda (t-s)} \left\|\phi_{\alpha}(s)  -\phi_{\beta}(s) \right\|ds \leq 
\displaystyle \frac{1}{\lambda} \sup_{t \leq \theta_{k_0}} \left\|\phi_{\alpha}(t) -\phi_{\beta}(t)\right\| 
\end{eqnarray*}
we obtain that
\begin{eqnarray*}
\left\| \phi_{\alpha}(t)-\phi_{\beta}(t)\right\|  < \displaystyle \frac{N(L_1+L_2)}{\lambda}  \sup_{t \leq \theta_{k_0}} \left\|\phi_{\alpha}(t) -\phi_{\beta}(t)\right\|  + \frac{N\eta\varepsilon}{\lambda}.
\end{eqnarray*}
Thus, the inequality
$$\sup_{t \leq \theta_{k_0}} \left\|\phi_{\alpha}(t) -\phi_{\beta}(t)\right\| \leq \displaystyle \frac{N \eta \varepsilon}{\lambda - N(L_1+L_2)} <\varepsilon $$
is valid. Accordingly, $\left\|\phi_{\alpha}(t) -\phi_{\beta}(t)\right\| \to 0$ as $t\to -\infty$.
Hence, $\phi_{\beta}(t) \in W^u(\phi_{\alpha} (t))$. $\square$

The main results of the present paper is mentioned in the subsequent theorem.

\begin{theorem} \label{mainthm}
Under the assumptions $(A1)-(A3)$ and $(A5)$, the following assertions are fulfilled: 
\begin{itemize}
	\item [i.] If a sequence $\beta \in \Lambda$ is homoclinic to a sequence $\alpha \in \Lambda$, then the function $\phi_{\beta} (t)\in \Omega$ is homoclinic to the function $\phi_{\alpha} (t)\in \Omega$.
	\item[ii.] If a sequence $\beta \in \Lambda$ is heteroclinic to  sequences $\alpha^1, \alpha^2 \in \Lambda$, then the function $\phi_{\beta} (t)\in \Omega$ is heteroclinic to the functions $\phi_{\alpha^1}(t), \phi_{\alpha^2} (t)\in \Omega$.
	\item[iii.] If the set $\Lambda$ is hyperbolic, then the set $\Omega$ is also hyperbolic.
\end{itemize}
\end{theorem}

\noindent \textbf{Proof.}
\begin{itemize}
\item [i.] Since $\beta$ is homoclinic to $\alpha$, we have $\beta \in W^s(\alpha) \cap W^u(\alpha)$. Then, $\beta$ belongs to both $W^s(\alpha)$ and $W^u(\alpha)$. According to Lemma \ref{lemmahomoclinic1}, $\phi_{\beta}(t)$ takes place in the set $W^s\left( \phi_{\alpha}(t)\right)$. On the other hand, Lemma \ref{lemmahomoclinic2} implies that $\phi_{\beta}(t)$  is an element of $W^u\left( \phi_{\alpha}(t)\right)$. Thus, $\phi_{\beta}(t) \in W^s\left( \phi_{\alpha}(t)\right)\cap W^u\left( \phi_{\alpha}(t)\right)$. Consequently, $\phi_{\beta}(t)$ is homoclinic to $\phi_{\alpha}(t)$.
\item [ii.] Being heteroclinic to $\alpha^1$, $\alpha^2$, the sequence $\beta$ lies in both of the sets $W^s(\alpha^1)$ and $W^u(\alpha^2)$. It follows from Lemma \ref{lemmahomoclinic1} that $\phi_{\beta}(t) \in W^s(\phi_{\alpha^1}(t))$. Additionally, $\phi_{\beta}(t) \in W^u(\phi_{\alpha^2}(t))$ by Lemma \ref{lemmahomoclinic2}. For that reason, $\phi_{\beta}(t) \in W^s(\phi_{\alpha^1}(t)) \cap W^u(\phi_{\alpha^2}(t))$. Consequently, $\phi_{\beta} (t)$ is heteroclinic to  $\phi_{\alpha^1}(t), \phi_{\alpha^2} (t)$.
\item [iii.] Fix an arbitrary sequence $\alpha \in \Lambda$. Due to the hyperbolicity of $\Lambda$, there exist sequences $\beta^1, \beta^2 \in \Lambda$, both of which are different from $\alpha$, such that $\beta^1 \in W^s(\alpha)$ and $\beta^2 \in W^u(\alpha)$.  Lemma \ref{lemmahomoclinic1} and Lemma \ref{lemmahomoclinic2} respectively imply that $\phi_{\beta^1}(t) \in W^s(\phi_{\alpha}(t))$ and $\phi_{\beta^2}(t) \in W^u(\phi_{\alpha}(t))$. Thus, $\Omega$ is hyperbolic. $\square$
\end{itemize}

An example which supports Theorem \ref{mainthm} is provided in the next section, and it is based on the logistic map.

\section{An Example} \label{sec4}

Let us consider the system
\begin{eqnarray} \label{systemexample}
&& z'_1(t) = 2z_1(t) -2z_2(t) + 0.03\cos(z_1(t))- 0.01\sin(z_2(\gamma(t))) + \frac{e^t}{1+e^t} + g_0(t,p^1) \nonumber \\
&& z'_2(t) = 5z_1(t)-3z_2(t) + 0.02\sin(z_2(t))+ 0.01\cos(z_1(\gamma(t))) +g_0(t,p^2),
\end{eqnarray}
in which $t\in\mathbb R$, $\theta_k=3k/2$ and $\zeta_k=(3k+1)/2$ for each $k\in\mathbb Z$, $\gamma(t)=\zeta_k$, $g_0(t,p^1) = p^1_k$,  $g_0(t,p^2) = p^2_k$ for $\theta_k \leq t <\theta_{k+1}$, $k\in\mathbb Z$, and the sequences $p^1=\left\lbrace p^1_k\right\rbrace_{k\in\mathbb Z}$, $p^2=\left\lbrace p^2_k\right\rbrace_{k\in\mathbb Z}$ are orbits of the logistic map
\begin{eqnarray} \label{systemdiscrete}
\eta_{k+1} = F_{\mu}(\eta_k), 
\end{eqnarray}
where $F_{\mu}(s)=\mu s(1-s)$, $s\in[0,1]$, and $\mu \in (0,4]$ is a parameter. The unit interval $[0,1]$ is invariant under the iterations of (\ref{systemdiscrete}) for those values of $\mu$ \cite{Hale91}.

System (\ref{systemexample}) is in the form of (\ref{mainsystem}) with 
$$z(t) = \begin{pmatrix} z_1(t) \\ z_2(t) \end{pmatrix}, \ A= \begin{pmatrix} 2 && -2 \\ 5 && -3 \end{pmatrix}, $$
$$f(t,z(t), z(\gamma(t))) = \begin{pmatrix} 0.03\cos(z_1(t))- 0.01\sin(z_2(\gamma(t))) + \displaystyle \frac{e^t}{1+e^t} \\ 0.02\sin(z_2(t))+ 0.01\cos(z_1(\gamma(t))) \end{pmatrix},$$
$$g(t,\alpha) = \begin{pmatrix} g_0(t,p^1) \\ g_0(t,p^2) \end{pmatrix},$$ and
$$\alpha = \left\lbrace \alpha_k\right\rbrace_{k\in\mathbb Z},\  \alpha_k = \begin{pmatrix} p_k^1 \\ p_k^2 \end{pmatrix}.$$
The eigenvalues of the matrix $A$ are $-1/2 \pm i\sqrt {15}/2$, and 
$e^{At}=Pe^{Bt}P^{-1}$, where
$$P= \begin{pmatrix} 0 && 4 \\ -\sqrt{15} && 5 \end{pmatrix}, \ B= \begin{pmatrix} -1/2 && -\sqrt{15}/2 \\ \sqrt{15}/2 && -1/2 \end{pmatrix}.$$
The inequality $\left\|e^{At} \right\|\leq N e^{-\lambda t}$ holds for $t\geq 0$, where $N=\left\| P\right\| \left\| P^{-1}\right\| = \left(7 + \sqrt{34} \, \right)/\sqrt{15}$ and $\lambda=1/2$.
The assumptions $(A1)-(A5)$ are fulfilled for system (\ref{systemexample}) with $\omega=3/2$, $M_f = 1.07229$, $L_1 =0.03$, and $L_2=0.01$. 

For a fixed value of the parameter $\mu \in (0,4]$, on the intervals $[0,1/2]$ and $[1/2,1]$ the function $F_{\mu}(s)$ admits  the inverses $$G_{\mu}(s) = \left(1-\sqrt{1-4s/\mu} \right)/2$$ and
$$H_{\mu}(s) = \left(1+\sqrt{1-4s/\mu} \right)/2,$$ respectively.

It was demonstrated by Avrutin et al. \cite{Avrutin15} that for $\mu=3.9$ the orbit 
$$\left\lbrace \ldots, H^3_{3.9}(1/3.9), H^2_{3.9}(1/3.9), H_{3.9}(1/3.9), 1/3.9, F_{3.9}(1/3.9), F^2_{3.9}(1/3.9), F^3_{3.9}(1/3.9),\ldots \right\rbrace$$
of the logistic map (\ref{systemdiscrete}) is homoclinic to its fixed point $2.9/3.9$. For that reason the sequence $\beta=\left\lbrace \beta_k\right\rbrace_{k\in\mathbb Z}$ is homoclinic to $\alpha=\left\lbrace \alpha_k\right\rbrace_{k\in\mathbb Z}$, where $$\alpha_k =  \begin{pmatrix} 2.9/3.9 \\ 2.9/3.9 \end{pmatrix}, \, k\in\mathbb Z,$$  
$$\beta_k = \begin{pmatrix} F^k_{3.9}(1/3.9) \\ F^k_{3.9}(1/3.9) \end{pmatrix}, \, k\geq 0, $$  and
$$\beta_k = \begin{pmatrix} H^{-k}_{3.9}(1/3.9) \\ H^{-k}_{3.9}(1/3.9) \end{pmatrix}, \, k< 0.$$ 
Therefore, according to Theorem  \ref{mainthm}, the bounded solution $\phi_{\beta}(t)$ of system (\ref{systemexample}) is homoclinic to the bounded solution $\phi_{\alpha}(t)$.
 
Next, let us investigate the presence of an heteroclinic solution of (\ref{systemexample}) by taking $\mu=4$. The orbit
$$\left\lbrace \ldots, G^3_{4}(1/4), G^2_{4}(1/4), G_{4}(1/4), 1/4, F_{4}(1/4), F^2_{4}(1/4), F^3_{4}(1/4),\ldots \right\rbrace$$
of (\ref{systemdiscrete}) is heteroclinic to the fixed points $3/4$ and $0$ \cite{Avrutin15}. Accordingly, the sequence $\widetilde{\beta} = \left\lbrace \widetilde{\beta}_k \right\rbrace _{k\in\mathbb Z}$ is heteroclinic to $\alpha^1 = \left\lbrace \alpha^1_k \right\rbrace _{k\in\mathbb Z}$, $\alpha^2 = \left\lbrace \alpha^2_k \right\rbrace _{k\in\mathbb Z}$, where
$$\alpha^1_k =  \begin{pmatrix} 3/4 \\ 3/4 \end{pmatrix}, \, \alpha^2_k =  \begin{pmatrix} 0 \\ 0 \end{pmatrix}, \, k\in\mathbb Z,$$ 
$$\widetilde{\beta}_k = \begin{pmatrix} F^k_{4}(1/4) \\ F^k_{4}(1/4) \end{pmatrix}, \, k\geq 0, $$  and
$$\widetilde{\beta}_k = \begin{pmatrix} G^{-k}_{4}(1/4) \\ G^{-k}_{4}(1/4) \end{pmatrix}, \, k< 0.$$ 
The bounded solution $\phi_{\widetilde{\beta}}(t)$ of system (\ref{systemexample}) is heteroclinic to the bounded solutions $\phi_{\alpha^1}(t)$, $\phi_{\alpha^2}(t)$ by Theorem  \ref{mainthm}.

\section{Conclusion} \label{sec5}

Long term behavior of solutions is one of the interesting subjects in the theory of differential equations. The present study makes a contribution to the theory in this sense such that a class of EPCAG is investigated from the asymptotic point of view. More precisely, the existence of homoclinic and heteroclinic solutions is demonstrated under sufficient conditions. The right-hand side of our model is discontinuous and is controlled by means of a discrete-time equation. The obtained results reveal that discontinuous systems with piecewise constant arguments are capable of generating homoclinic and heteroclinic motions. It is also deduced that the hyperbolic structure of a discrete map can be transmitted to a continuous-time dynamics through a system of type (\ref{mainsystem}). In the future our results can be developed for partial differential equations with piecewise constant arguments \cite{Stefanidou23,Veloz15}.

\section*{Appendix}

A proof of Lemma \ref{bddsolnlemma} is as follows.

\noindent \textbf{Proof of Lemma \ref{bddsolnlemma}.} Suppose that $\mathcal{C}$ is the set of all uniformly bounded continuous functions $\psi:\mathbb R \to \mathbb R^m$ such that $\left\|\psi \right\|_{\infty} \leq M_{\phi}$, where $\left\|\psi \right\|_{\infty} = \displaystyle \sup_{t \in \mathbb R} \left\|\psi(t)\right\|$ and $M_{\phi}$ is the number defined by (\ref{numbermphi}).

Let a sequence $\alpha = \left\lbrace \alpha_k \right\rbrace_{k\in\mathbb Z} \in \Lambda$ be given. We consider the operator $\Pi$ defined on $\mathcal{C}$ by the equation
$$
(\Pi(\psi))(t) = \displaystyle \int_{-\infty}^t e^{A(t-s)} \left[ f(s,\psi(s), \psi(\gamma (s))) + g(s,\alpha)\right] ds.
$$
For a given function $\psi(t) \in \mathcal{C}$, we have for every $t \in \mathbb R$ that
\begin{eqnarray*}
\left\| (\Pi(\psi))(t)\right\| & \leq & \displaystyle \int_{-\infty}^t Ne^{-\lambda(t-s)}  \left\| f(s,\psi(s), \psi(\gamma (s))) + g(s,\alpha)  \right\|  ds \\
& \leq & N(M_f + M_F)  \displaystyle \int_{-\infty}^t  e^{-\lambda(t-s)}ds \\ &=& M_{\phi}.
\end{eqnarray*}
The last inequality yields $\left\|  \Pi(\psi)  \right\|_{\infty} \leq \left\| \psi\right\|_{\infty} $. Hence, $  \Pi \left(  \mathcal{C}\right)  \subseteq  \mathcal{C}$.

In addition to that, if $\psi_1(t)$ and $\psi_2(t)$ are elements of $\mathcal{C}$, then using the equation
\begin{eqnarray*}
(\Pi(\psi_1))(t) - (\Pi(\psi_2))(t) = \displaystyle \int_{-\infty}^t e^{A(t-s)} \left[ f(s,\psi_1(s), \psi_1(\gamma (s)))  - f(s,\psi_2(s), \psi_2(\gamma (s))) \right] ds
\end{eqnarray*}
one can attain that
\begin{eqnarray*}
\left\| (\Pi(\psi_1))(t) - (\Pi(\psi_2))(t)\right\| & \leq & N L_1  \displaystyle \int_{-\infty}^t e^{-\lambda (t-s)} \left\| \psi_1(s) - \psi_2(s)\right\| ds \\
&+&  N L_2  \displaystyle \int_{-\infty}^t e^{-\lambda (t-s)} \left\| \psi_1( \gamma (s)) - \psi_2( \gamma (s) )\right\| ds \\
& \leq & \displaystyle \frac{N(L_1+L_2)}{\lambda} \left\|\psi_1-\psi_2 \right\|_{\infty}. 
\end{eqnarray*}
Therefore, the inequality
$$\left\| (\Pi(\psi_1)) - (\Pi(\psi_2))\right\|_{\infty} \leq  \displaystyle \frac{N(L_1+L_2)}{\lambda} \left\|\psi_1-\psi_2 \right\|_{\infty}$$
is valid. The operator $\Pi$ is a contraction in accordance with the assumption $(A4)$. Thus, for every sequence $\alpha = \left\lbrace \alpha_k \right\rbrace_{k\in\mathbb Z} \in \Lambda$, there exists a unique bounded solution $\phi_{\alpha}(t)$ of system 
(\ref{mainsystem}) satisfying $\displaystyle \sup_{t\in\mathbb R} \left\| \phi_{\alpha}(t)\right\| \leq M_{\phi}$. $\square$

%
%
%
%
%
%
%
%
%
%
%
%
%


\end{document}